\documentclass[lettersize,journal]{IEEEtran}
\usepackage{amsmath,amsfonts}
\usepackage{algorithmic}
\usepackage{algorithm}
\usepackage{array}
\usepackage[caption=false,font=normalsize,labelfont=sf,textfont=sf]{subfig}
\usepackage{textcomp}
\usepackage{url}
\usepackage{verbatim}
\usepackage{graphicx}
\usepackage{cite}
\usepackage{hyperref}
\hypersetup{
	colorlinks=true, 
	linkcolor=blue,  
	urlcolor=blue,   
	citecolor=blue   
}

\usepackage{xcolor,soul,framed} 
\colorlet{shadecolor}{yellow}
\usepackage{eqparbox}
\usepackage{url}

\usepackage{multicol}
\usepackage{listings}
\usepackage{answers}
\usepackage{setspace}
\usepackage{enumitem}
\usepackage{multicol}
\usepackage{mathrsfs}
\usepackage{amsmath,amsthm,amssymb}
\usepackage{esint}
\usepackage{textpos}
\usepackage{dblfloatfix}
\usepackage{xcolor}
\usepackage{xintbinhex}
\usepackage{calc}
\usepackage{mathtools}
\usepackage{epstopdf}
\usepackage{comment}

\newcommand{\R}{\mathbb{R}}

\newcommand{\vn}[1]{\ensuremath{\boldsymbol{#1}}}                   
\newcommand{\vet}[1]{{{ \mathbf{#1}}}}                      
\newcommand{\mat}[1]{{{{\mathbf{{#1}}}}}}         
 
\newcommand{\diver}{\nabla\cdot}
\newcommand{\diverp}{\nabla'\cdot}

\makeindex             
\begin{document}

\title{A Modified Dielectric Contrast based Integral Equation for 2D TE Scattering by Inhomogeneous Domains}





\author{Akshay Pratap Singh,~\IEEEmembership{Student Member,~IEEE,}
	Kuldeep Singh,~\IEEEmembership{Member,~IEEE,}\\
	Rajendra Mitharwal,~\IEEEmembership{Member,~IEEE}
\thanks{Manuscript received July 19, 2025. This work is supported by the Institute Research Assistantship Grant at MNIT Jaipur, funded by the Ministry of Education, Government of India. \textit{(Corresponding Author: Akshay Pratap Singh)}}

\thanks{Akshay Pratap Singh, Kuldeep Singh, and Rajendra Mitharwal are	with Department of ECE, Malaviya National Insitute of Technology Jaipur, JLN Marg, Jaipur - 302017, India (e-mail: 2020rec9501@mnit.ac.in;	kuldeep.ece@mnit.ac.in; rajendra.ece@mnit.ac.in)}}

\IEEEpubid{0000--0000/00\$00.00~\copyright~2025 IEEE}

\maketitle

\begin{abstract}
This work presents a modified domain integral equation approach for the forward problem of TE scattering, employing a modified definition of dielectric contrast and discretizing the electric field density using Rao-Wilton-Glisson (RWG) basis functions. The proposed formulation mitigates the numerical challenges introduced by the gradient-divergence operator in traditional electric field-based vector formulations. The use of RWG basis functions over triangular meshes enhances geometric conformity, ensures tangential continuity across dielectric interfaces, and facilitates the application of well known singularity extraction techniques for numerical accuracy. Validation through numerical experiments on a two-layered dielectric cylinder demonstrates excellent agreement between computed and analytical scattered fields. Convergence studies confirm improving solution accuracy with mesh refinement indicating robustness with respect to discretization without increasing the iterations.
\end{abstract}

\begin{IEEEkeywords}
Domain Integral Equations, Greens function, Inhomogeneous dielectric media, RWG basis functions, TE Scattering.
\end{IEEEkeywords}

\section{Introduction}
\IEEEPARstart{M}{icrowaves} have been used prominently in characterizing electromagnetic materials using non destructive techniques especially for imaging civil structures and monitoring processes in industrial settings \cite{9000635,7803599, Catapano2007Simple, Crocco2012Linear}. To facilitate such characterization, a penetrable scatterer is illuminated by a known source typically Transverse Electric (TE) polarized waves, followed by the measurement of the resulting scattered fields \cite{Kooij1998Nonlinear}. The quality of imaging can be significantly improved by incorporating additional nonredundant scattered field data through an increased number of measurement antennas. However, the practical implementation of a larger antenna array is constrained by spatial limitations, substantial experimental costs, and extended data acquisition durations. A frequency-domain zero-padding (FDZP) interpolation technique can be used to facilitate the acquisition of scattered field data at additional virtual antenna positions \cite{10505037}. The dielectric contrast ratio enables the estimation of the smallest detectable target contrast for a given size or vice versa \cite{6889125}. These measured fields are then utilized as input to inverse scattering algorithms. These algorithms aim to reconstruct the material properties by formulating the relationship between the unknown parameters and the scattered fields as a two-step optimization problem: forward and inverse problem. The forward problem deals with computing scattered field due to the inhomogeneous dielectric scatterer illuminated by TE polarized wave. 

Mathematically, the forward problem can be formulated into a scalar domain integral equations in terms of unknown mangetic field component or a vector domain integral equation in terms of electric field vector (two unknown vector) \cite{lori2017new, franceschini2006iterative}. It has been demonstrated that the scalar integral equation using single magnetic field component has lower accuracy compared to vector domain integral equation based on two component electric field unknowns \cite{Kooij1998Nonlinear}. Moreover, the magnetic field–based integral equation tends to suffer from reduced numerical stability and efficiency \cite{mojabi2010comparison}. In contrast, although the electric field–based vector integral equation is generally more stable and accurate, it presents computational challenges due to the presence of the gradient-divergence operator\cite{mojabi2010comparison}. The prominent methods \cite{Kooij1998Nonlinear} address the singularity of the gradient-divergence operator by applying the involved spatial derivatives on the spatially continuous function.  This strategy ensures that the resulting numerical scheme avoids the generation of non-physical surface currents. This is achieved by employing roof-top basis functions \cite{gurel1997quantitative} in the direction of spatial derivatives. 

\IEEEpubidadjcol

In this work, we propose an alternative discretization strategy  wherein a modified definition of the dielectric contrast ratio is utilized compared to conventional definition \cite{Kooij1998Nonlinear}. This redefinition transforms the unknown quantity of the integral equation into a function of electric field density which is discretized using the Rao-Wilton-Glisson (RWG) vector basis functions \cite{Rao1982Electromagnetic}. These basis functions are naturally suited to triangular mesh discretizations, representing complex scatterer geometries more accurately. They inherently conform to the divergence operator and maintain the tangential continuity of electric field density across inhomogeneous dielectric interfaces.  These advantages are further enhanced by applying singularity extraction techniques\cite{wilton1984potential}, which improve the accuracy of computing the diagonal and near-interaction elements of the system matrix \cite{bao2020fully}. In addition, the use of RWG basis functions facilitates the integration of fast solution algorithms \cite{andriulli2012helmholtz} in the forward step of inverse scattering algorithms.

The following section introduces the necessary mathematical background for the development of the reformulated integral equation. The third section introduces the proposed reformulated integral equations and its discretization using RWG basis functions to handle the TE scattering problem in inhomogeneous dielectric media. The numerical results section outlines the numerical experiments conducted to validate the accuracy and effectiveness of the proposed formulation, followed by conclusion.

\section{Theoretical Background}
Consider an inhomogeneous dielectric region $\Omega \subset \R^2$, bounded by the curve $\Gamma$, and illuminated by a time-harmonic TE polarized plane wave with incident electric field $\vn{E}^{inc}(\vn{r}) \in \R^2$. The interaction with the dielectric medium gives rise to a scattered field $\vn{E}^s(\vn{r})$, such that the total field in the exterior of $\Omega$ is expressed as $\vn{E}(\vn{r}) = \vn{E}^i(\vn{r}) + \vn{E}^s(\vn{r})$.

This interaction can be obtained analytically \cite{balanis2024balanis} for a multi-layered cylinder in cylindrical coordinate system for the scattered field $\vn{E}^s(\vn{r}) = E^s_\rho(\rho,\phi) \hat{\rho} + E^s_\phi(\rho,\phi)\hat{\phi}$ in the external region of $\Omega$ as

\begin{equation}\label{eq:analytical}
	\begin{aligned}
		E_\rho^s(\rho,\phi)=&\frac{1}{j\omega\epsilon_0} \frac{1}{\rho} H_o \sum_{n=-\infty}^{\infty} a_n H_n^{(2)}(k_0 \rho)jn\exp(jn\phi)\\
		E_\phi^s(\rho,\phi)=&\frac{-\beta_0}{j\omega\epsilon_0} H_0 \sum_{n=-\infty}^{\infty} a_n H_n^{'(2)}(k_0 \rho)\exp(jn\phi)
	\end{aligned}
\end{equation}

where $\omega$ is the angular frequency, $\epsilon_0$ is the free space permittivity, $k_0$ is the free space wave number, $a_n$ are the expansion coefficients depending on the number of layers of the dielectric region in the cylinder, $H_0$ is the amplitude of the magnetic field of the incident plane wave, and $H_n^{(2)}$ represents Hankel function of second kind of order $n$ \cite{zhang1997computation}.

\begin{figure}
	\centering
	\includegraphics[width=0.25\textwidth] {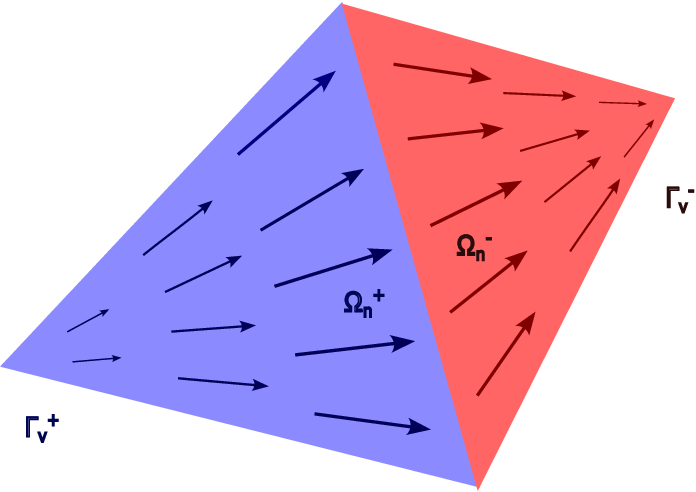}
	\caption{An RWG basis function defined on pair of adjacent triangles}
	\label{fig:rwg}
\end{figure}

However, when the $\Omega$ represents a complex geometry, the interaction needs to be formulated mathematically using the integral equation\cite{Kooij1998Nonlinear} as

\begin{equation}\label{eq:vandenberg}
	\vn{E}^{inc}(\vn{r}) = \vn{E}(\vn{r})-(k_0^2+\nabla\nabla)\int_{\Omega}\chi(\vn{r}')\vn{E}(\vn{r}')G(\vn{r},\vn{r}')d\Omega', 
\end{equation}
where $\chi(\vn{r}) = \frac{\epsilon(\vn{r})-\epsilon_0}{\epsilon_0}$ is the dielectric contrast defined in terms of the spatially varying electric permittivity  $\epsilon(\vn{r})$ within $\Omega$ and the free-space permittivity $\epsilon_0$, and  $G(\vn{r},\vn{r}') = \frac{j}{4}H_0^{(2)}(k_0|\vn{r}-\vn{r}'|)$ being the 2D Green's function \cite{gibson2021method}. The integral term in the above equation represents the scattered electric field resulting from the dielectric region $\Omega$.
In the \cite{Kooij1998Nonlinear}, the electric field $\vn{E}(\vn{r}')$, which is the unknown in the equation, is discretized using a pair of two-dimensional scalar basis functions (along the $x$ and $y$ directions), typically defined using pulse and rooftop functions on a rectangular mesh. In this work, we propose a modified formulation wherein the unknown function in the integral equation is discretized using RWG vector basis functions on a triangular mesh. Each edge shared by a pair of adjacent triangles in this mesh defines an RWG basis function, as shown in Figure \ref{fig:rwg}. The triangles $\Omega_n^+$ and $\Omega_n^-$ represent the positive and negative terminals, respectively, with the vector basis function oriented from $\Omega_n^+$ to $\Omega_n^-$.  Mathematically, the RWG basis function is expressed as
\begin{equation}\label{eq:rwg}
	\vn{f}_n (\vn{r})= 
	\begin{dcases}						
		\frac{\vn{r}-\vn{r_{v}^+}}{2A_{n}^+}, & \vn{r} \in \Omega_n^+\\
		-\frac{\vn{r}-\vn{r_{v}}^-}{2A_{n}^-},&\vn{r} \in \Omega_n^-
	\end{dcases}
\end{equation}
defined on $n^{\text{th}}$ edge with $\vn{r_{v}^+}$ and $\vn{r_{v}^-}$ denoting the reference vertices of the triangles $\Omega_n^+$ and $\Omega_n^-$, respectively, and $A_{n}^+$ and $A_{n}^-$ represent their corresponding areas, as shown in Figure~\ref{fig:rwg}.

\section{Proposed Reformulation}
The \eqref{eq:vandenberg} can be reformulated by adopting an alternative definition of dielectric contrast ratio 
\begin{equation}
	\chi(\vn{r}) = \frac{\epsilon(\vn{r})-\epsilon_0}{\epsilon_0}
\end{equation}
(different from \cite{Kooij1998Nonlinear} ) as
\begin{equation*}\label{23}
	\begin{aligned}
		\vn{E}^{inc}(\vn{r}) = & \vn{E}(\vn{r})+\frac{jk_0}{\epsilon_0}(j k_0-\frac{\nabla\nabla}{j _0})\int_{\Omega} G(\vn{r},\vn{r}')\chi(\vn{r}')\\
		& \qquad \epsilon(\vn{r}')\vn{E}(\vn{r}')dr',
	\end{aligned}
\end{equation*}
which can be simplified in terms of following electric current density in terms of electric flux density $\vn{D}(\vn{r}) = \epsilon(\vn{r})\vn{E}(\vn{r})$ 
\begin{equation}\label{eq:currentfieldrel}
	\vn{J}(\vn{r}) = j\omega\eta_0 \vn{D}(\vn{r}) = \frac{jk_0}{1 - \chi (\vn{r})}\vn{E} (\vn{r}),
\end{equation}
($\eta_0$ being the wave impedance of free space) as follows

\begin{equation}\label{eq:sie}
	\begin{aligned}
		\vn{E}^{inc}(\vn{r}) &= \frac{\vn{J}(\vn{r}) }{jk_0}(1-\chi(\vn{r}))\\
		&\qquad +{jk_0}\int_{\Omega} G(\vn{r},\vn{r}') {\chi(\vn{r}')}\vn{J}(\vn{r}') d\Omega'	\\
		&\qquad\qquad-\frac{\nabla\nabla}{jk_0}\int_{\Omega} G(\vn{r},\vn{r}') {\chi(\vn{r}')}\vn{J}(\vn{r}') d\Omega'.
	\end{aligned}
\end{equation}
The permittivity $\epsilon(\vn{r})=\epsilon_r(\vn{r})\epsilon_0-j\frac{\sigma(\vn{r})}{\omega}$ used in the dielectric contrast ratio can be a complex quantity due to conductivity $\sigma(\vn{r})$ and space varying relative permittivity $\epsilon_r(\vn{r})$.

\subsection {Discretization of the reformulated integral equation}
The scatterer domain $\Omega$ is discretized using a triangular mesh, where a RWG vector basis function is defined on each interior edge shared by a pair of adjacent triangles. The unknown surface current density in the domain integral equation \eqref{eq:sie} is discretized using RWG vector basis functions (defined in \eqref{eq:rwg}) which provide an appropriate representation for piecewise linear currents over triangular meshes, as follows
\begin{equation}\label{eq:rwgdiscretization}
	\vn{J}(\vn{r})=\sum_{n=1}^{N}d_{n} \vn{f}_n (\vn{r}),
\end{equation}		
where $N$ denotes the total number of interior edges within the domain $\Omega$. Substituting the discretized surface current into \eqref{eq:sie} and applying the Galerkin method yields the following discretized linear system of equations:
\begin{equation}\label{eq:sysmat}
	(\mat{I_\chi} + \mat{Z_A} + \mat{Z_\phi} ) \vet{d} = \vet{e}.
\end{equation}
The element-wise definitions for $m^{\text{th}}$ row and $n^{\text{th}}$ column of the involved matrices are given by
\begin{equation}
	\mat{(I_\chi)}_{mn} = \frac{1}{jk_0} \int_{\Omega_m} \vn{f}_m(\vn{r})\cdot \vn{f}_n(\vn{r})  (1-\chi(\vn{r})) d\Omega,
\end{equation}
\begin{equation}
	\mat{(Z_A)}_{mn} = {jk_0}\int_{\Omega_m}\vn{f}_m(\vn{r})\cdot\int_{\Omega_n} G(\vn{r},\vn{r}') {\chi(\vn{r}')}\vn{f}_n(\vn{r'}) d\Omega'd\Omega,
\end{equation}
and
\begin{multline}
		\mat{(Z_\phi)}_{mn} = \frac{1}{jk_0}\int_{\Omega_m}\diver\vn{f}_m(\vn{r})\int_{\Omega_n} G(\vn{r},\vn{r}') \\  \diverp\{{\chi(\vn{r}')}\vn{f}_n(\vn{r'})\} d\Omega'd\Omega
\end{multline}
where $\Omega_m$ denotes the support of $m^{\text{th}}$ RWG basis function $\vn{f}_m(\vn{r})$ on testing side with corresponding interpretation for the $n^{\text{th}}$ RWG basis function on source side. The last expression represents the scalar potential contribution incoprorating the divergence term in the source integral over $\Omega_n$ which can be further decomposed into two terms \cite{kobidze2004integral} using
\begin{equation}\label{eq:divchirwg}
	\diver\{{\chi(\vn{r})}\vn{f}_n(\vn{r})\}= 
	\begin{dcases}						
		\frac{\chi^{\pm}_n}{A_{n}^\pm}, & \vn{r} \in \Omega_n^\pm\\
		-\delta(\vn{r})\chi^{\pm} \hat{n}^{\pm}\cdot \vn{f}_n(\vn{r}),&\vn{r} \in \Gamma_n
	\end{dcases}
\end{equation}
where $\hat{n}^{\pm}$ are outward directed normal of the reference triangles (with a common edge $\Gamma_n$) of the $n^{\text{th}}$ RWG basis function, and $\delta(\vn{r})$ is the delta functional whose support is $\Gamma_n$. The sparse matrix $\mat{(I_\chi)}$ can be computed analytically by iterating over the set of triangles discretizing $\Omega$. For example, consider the $p^{\text{th}}$ triangle $\Omega_p$ (defined by the vertex coordinates $\vn{v}_1, \vn{v}_2, \text{ and } \vn{v_3}$ having an area of $A_p$) with each of its edges associated with an RWG basis function. Assume $\vn{f}_m(\vn{r})$ to be an RWG function that has an outward orientation $s_m$  over the edge of $\Omega_p$ with a corresponding value $+1$ or $-1$. Using $\vn{v}_{31} = \vn{v}_3 - \vn{v}_1, \text{ and } \vn{v}_{21} = \vn{v}_2 - \vn{v}_1$, $\Omega_p$ will contribute upto nine elements update as $\mat{(I_\chi)}_{mn}= \mat{(I_\chi)}_{mn} + s_m s_n\frac{(1-\chi(\vn{r}))}{jk_0} c_p$ where
\begin{equation}\label{eq:rwgc_p}
	c_p = 
	\begin{dcases}						
		\frac{1}{24 A_p} (|\vn{v}_{31}|^2 + |\vn{v}_{21} |^2+ \vn{v}_{31}\cdot \vn{v}_{21}), & m = n\\
			\frac{1}{24 A_p} (|\vn{v}_{31}|^2 + |\vn{v}_{21} |^2+ \vn{v}_{31}\cdot \vn{v}_{21})\\ 
			\qquad- \frac{1}{12 A_p} (\vn{v}_{31} + \vn{v}_{21})\cdot \vn{v}_{21}
		,& m\ne n.
	\end{dcases},
\end{equation}
The matrix elements corresponding to scalar and vector potential terms can be computed using the Gaussian quadrature rules defined for triangular domain \cite{andriulli2012helmholtz}. Specifically, assume $\{\vn{r}_i,w_i\}$ represent a 2D quadrature rule of $N_G$ points over a triangle  $\Omega_p$. Then, the matrix element corresponding to vector potential term is
\begin{multline}\label{eq:matvecpot}
	 \mat{(Z_A)}_{mn} =\\ {jk_0}\left\{\sum_{\substack{i=1\\\vn{r}_i\in\Omega_m^+}}^{N_G}w_i(\vn{r}_i-\vn{r_{m}^+})\cdot\sum_{\substack{j=1\\\vn{r}_j\in\Omega_n^+}}^{N_G} w_j{\chi^{+}_n}G(\vn{r}_i,\vn{r}_j ) (\vn{r}_j-\vn{r_{n}^+})\right.\\
	 -\sum_{\substack{i=1\\\vn{r}_i\in\Omega_m^+}}^{N_G}w_i(\vn{r}_i-\vn{r_{m}^+})\cdot\sum_{\substack{j=1\\\vn{r}_j\in\Omega_n^-}}^{N_G}w_j {\chi^{-}_n} G(\vn{r}_i,\vn{r}_j ) (\vn{r}_j-\vn{r_{n}^-})\\
	-\sum_{\substack{i=1\\\vn{r}_i\in\Omega_m^-}}^{N_G}w_i(\vn{r}_i-\vn{r_{m}^-})\cdot\sum_{\substack{j=1\\\vn{r}_j\in\Omega_n^+}}^{N_G}w_j {\chi^{+}_n}G(\vn{r}_i,\vn{r}_j ) (\vn{r}_j-\vn{r_{n}^+})\\
	\quad \left.+\sum_{\substack{i=1\\\vn{r}_i\in\Omega_m^-}}^{N_G}w_i(\vn{r}_i-\vn{r_{m}^-})\cdot\sum_{\substack{j=1\\\vn{r}_j\in\Omega_n^-}}^{N_G}w_j {\chi^{-}_n} G(\vn{r}_i,\vn{r}_j ) (\vn{r}_j-\vn{r_{n}^-})\right\}
\end{multline}
For the matrix corresponding to the scalar potential term, it is also necessary to define a 1D quadrature rule $\{\vn{r}_i, w_i\}$ with $N_l$ points along the edge $\Gamma_n$ on which the RWG basis function $\Omega_n$ is defined. Using this, the corresponding matrix is
\begin{equation}\label{eq:matscalpot}
	\begin{aligned}
	 \mat{(Z_\phi)}_{mn} &= 
	  {jk_0}\left\{\sum_{\substack{i=1\\\vn{r}_i\in\Omega_m^+}}^{N_G} \sum_{\substack{j=1\\\vn{r}_j\in\Omega_n^+}}^{N_G} w_iw_jG(\vn{r}_i,\vn{r}_j ) {\chi^{+}_n}\right.	\\
	 &\qquad -\sum_{\substack{i=1\\\vn{r}_i\in\Omega_m^+}}^{N_G}\cdot\sum_{\substack{j=1\\\vn{r}_j\in\Omega_n^-}}^{N_G} w_iw_jG(\vn{r}_i,\vn{r}_j ) {\chi^{-}_n}
	\\
	&\qquad -\sum_{\substack{i=1\\\vn{r}_i\in\Omega_m^-}}^{N_G}\sum_{\substack{j=1\\\vn{r}_j\in\Omega_n^+}}^{N_G} w_iw_jG(\vn{r}_i,\vn{r}_j ) {\chi^{+}_n}\\ 
	&\qquad +\sum_{\substack{i=1\\\vn{r}_i\in\Omega_m^-}}^{N_G}\sum_{\substack{j=1\\\vn{r}_j\in\Omega_n^-}}^{N_G} w_iw_jG(\vn{r}_i,\vn{r}_j ) {\chi^{-}_n} \\
	&\qquad -\sum_{\substack{i=1\\\vn{r}_i\in\Omega_m^+}}^{N_G}\sum_{\substack{j=1\\\vn{r}_j\in\Gamma_n}}^{N_l} w_iw_jG(\vn{r}_i,\vn{r}_j ) ({\chi^{+}_n} - {\chi^{-}_n})\\ 
	&\qquad \left.+\sum_{\substack{i=1\\\vn{r}_i\in\Omega_m^-}}^{N_G}\sum_{\substack{j=1\\\vn{r}_j\in\Gamma_n}}^{N_l} w_iw_jG(\vn{r}_i,\vn{r}_j ) ({\chi^{+}_n} - {\chi^{-}_n})
	\right\}.
	\end{aligned}
\end{equation}
 \eqref{eq:matvecpot} and \eqref{eq:matscalpot} exhibit singular behavior when the quadrature points $\vn{r_i}$ and $\vn{r_j}$ are in close proximity. This issue is addressed by decomposing the Green's function into singular and non-singular components. The non-singular part is computed numerically using the quadrature scheme described earlier, while the singular component is treated analytically using the standard methods\cite{wilton1984potential}. The elements of right hand side vector is computed as
\begin{equation*}
		{\vn{e}_m} = \int_{\Omega_m}\vn{f}_m(\vn{r})\cdot\vn{E}_{inc}(\vn{r})d\Omega\\
\end{equation*}

\begin{multline}
		\qquad\qquad{\vn{e}_m}=\sum_{\substack{i=1\\\vn{r}_i\in\Omega_m^+}}^{N_G}w_i(\vn{r}_i-\vn{r_{m}^+})\cdot\vn{E}_{inc} (\vn{r}_i)\\
		-\sum_{\substack{i=1\\\vn{r}_i\in\Omega_m^-}}^{N_G}w_i(\vn{r}_i-\vn{r_{m}^-})\cdot\vn{E}_{inc} (\vn{r}_i)
\end{multline}
	
The system of equations in \eqref{eq:sysmat} can be solved using iterative solvers to obtain the unknown coefficients in the expansion given by  \eqref{eq:rwgdiscretization}. Once these coefficients are determined, the scattered electric field in the exterior of the domain $\Omega$ can be evaluated using
\begin{multline}
		\vn{E}_s(\vn{r}) = \sum_{n=1}^{N}d_{n} \sum_{\substack{i=1\\\vn{r}_i\in\Omega_n^{\pm}}}^{N_G}w_i A_n^{\pm} \left[{jk_0}G(\vn{r},\vn{r}_i) {\chi(\vn{r}_i)}\vn{f}_n(\vn{\vn{r}_i})\right.\\
		-\left. \frac{1}{jk_0} \nabla G(\vn{r},\vn{r}_i) \&\diverp\{{\chi(\vn{r}_i)}\vn{f}_n(\vn{r}_i)\} \right]
\end{multline}

for brevity a compact notation is used. 
The reader may readily extend the source integral term of \eqref{eq:matscalpot} to explicitly expand the divergence term in the above expression. The electric field inside $\Omega$ can be directly obtained using \eqref{eq:currentfieldrel} and \eqref{eq:rwgdiscretization}.

\section{Numerical Results}
\begin{figure}
	\centering
	\includegraphics[width=0.5\textwidth] {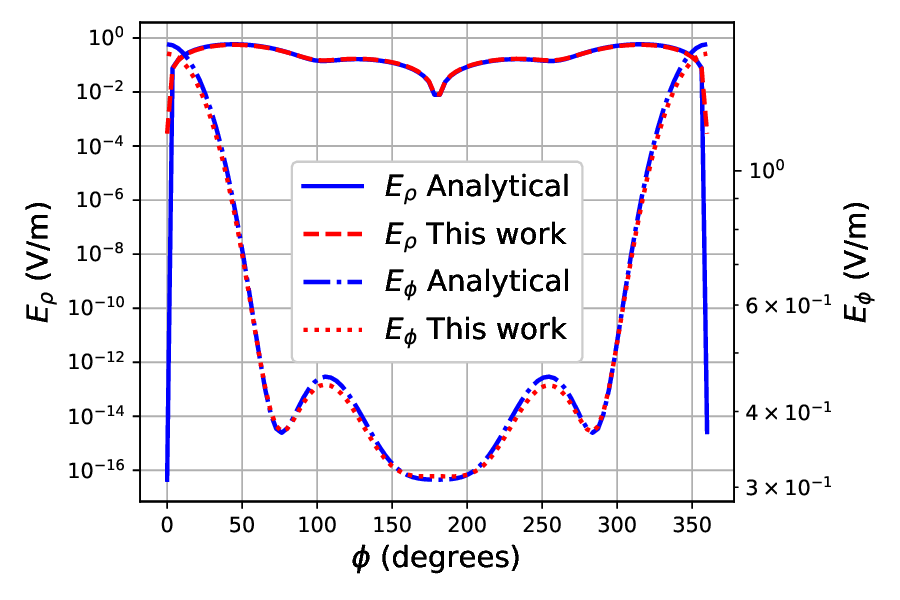}
	\caption{Comparison of the electric field components obtained using proposed method with the analytical expressions in cylindrical coordinates.}
	\label{fig:merged_Erho_Ephi_plot}
\end{figure}

\begin{figure}
	\centering
	\includegraphics[width=0.5\textwidth] {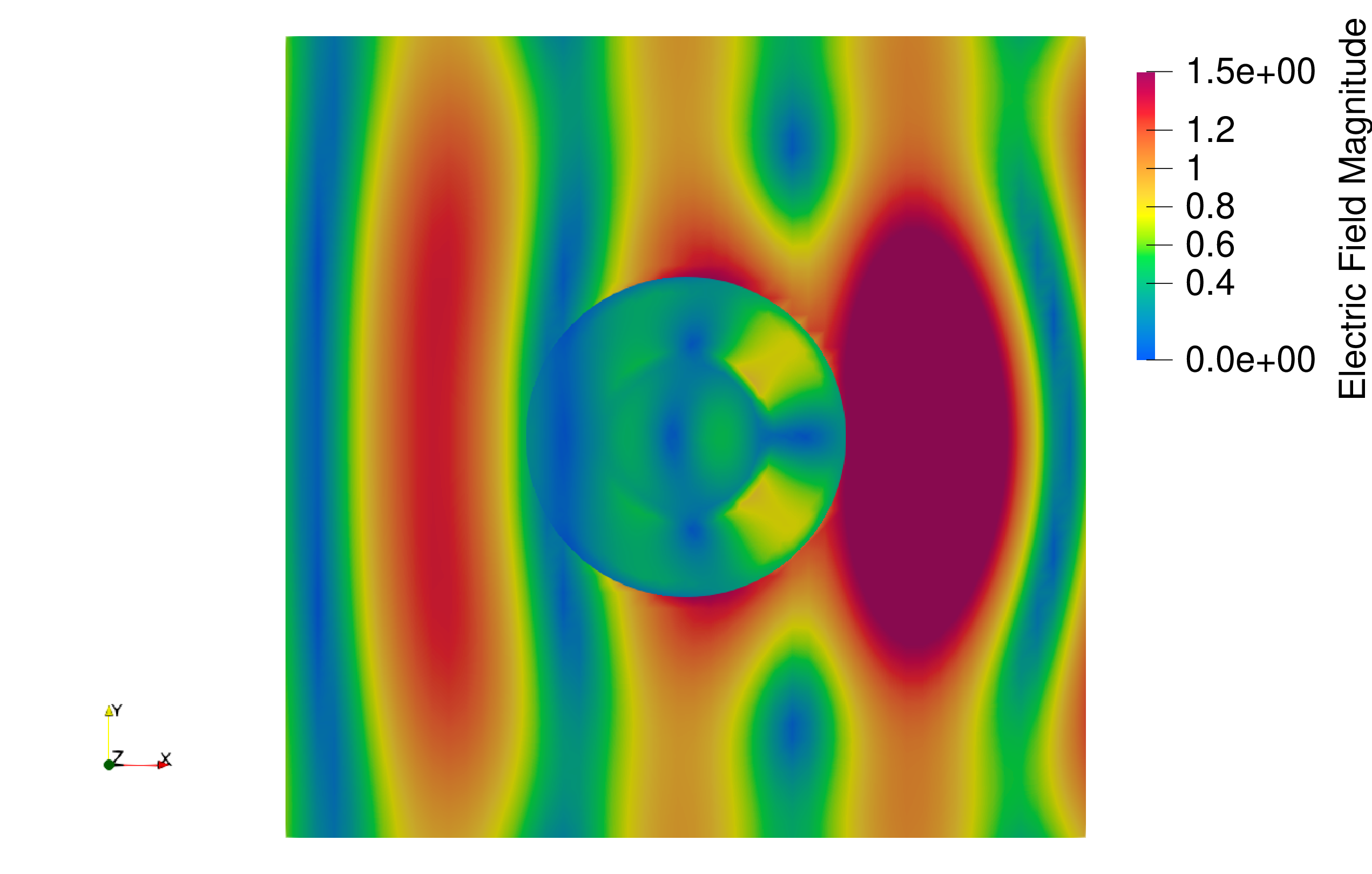}
	\caption{Electric field distribution inside and outside the two-layered dielectric cylinder at 1 GHz.}
	\label{fig:Total_electric_field_plot}
\end{figure}

The proposed reformulated integral equation is validated through its application to a two-layered dielectric cylinder in a two-dimensional setting. The cylinder has an inner radius of 0.05 meters and an outer radius of 0.1 meters. The relative permittivities of the inner and outer regions are set to 2 and 8, respectively. The structure is illuminated by a time-harmonic TE polarized plane wave at an operating frequency of 1 GHz, with the incident electric field given by $\vn{E}^{inc}(\vn{r}) = E_0 e^{ik_0x}\hat{x}$ with the amplitude $E_0=1$ V/m. The cylinder's domain $\Omega$ is discretized using a triangular mesh with an average edge length of 0.007 m, and numerical integration is performed using a single Gaussian quadrature point per triangle. This discretization yields a total of 4104 RWG basis functions. The scattered electric field is computed on a circle of radius 0.15 m in the $x-y$ plane. The cylindrical components of the scattered field, $(E_\rho, E_\phi)$ demonstrate excellent overlap with the analytically computed electric field (using \eqref{eq:analytical}) for the same example as depicted in Fig.~\ref{fig:merged_Erho_Ephi_plot}, and the corresponding field distribution is illustrated in Fig.~\ref{fig:Total_electric_field_plot}.

\begin{figure}
	\centering
	\includegraphics[width=0.5\textwidth] {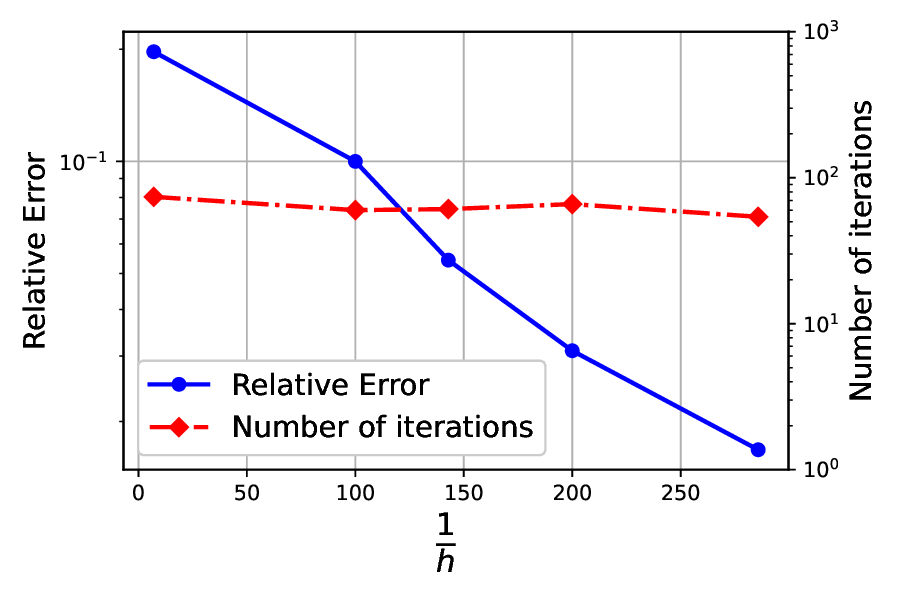}
	\caption{Relative error in the computed electric field and the number of solver iterations as functions of discretization density, where $h$ denotes the average edge length of the discretized geometrical mesh for the two layered dielectric cylinder.}
	\label{fig:h_convergence_plot}
\end{figure}

Next, the convergence behavior of the proposed solution is examined in Fig.~\ref{fig:h_convergence_plot} by refining the mesh by varying the average edge length from 0.141421 m to 0.0035 m. An improvement in the accuracy of the computed solution is observed as the mesh is refined. This is demonstrated by comparing the numerically obtained scattered field (evaluated on the same observation circle as in the previous example) with the corresponding analytical solution. The results confirm that an increase in the discretization density leads the numerically computed solution approximately more closer to the analytically obtained solution. The number of iterations required for the numerical solution remains nearly constant across different mesh densities.

\begin{figure}
	\centering
	\includegraphics[width=0.5\textwidth] {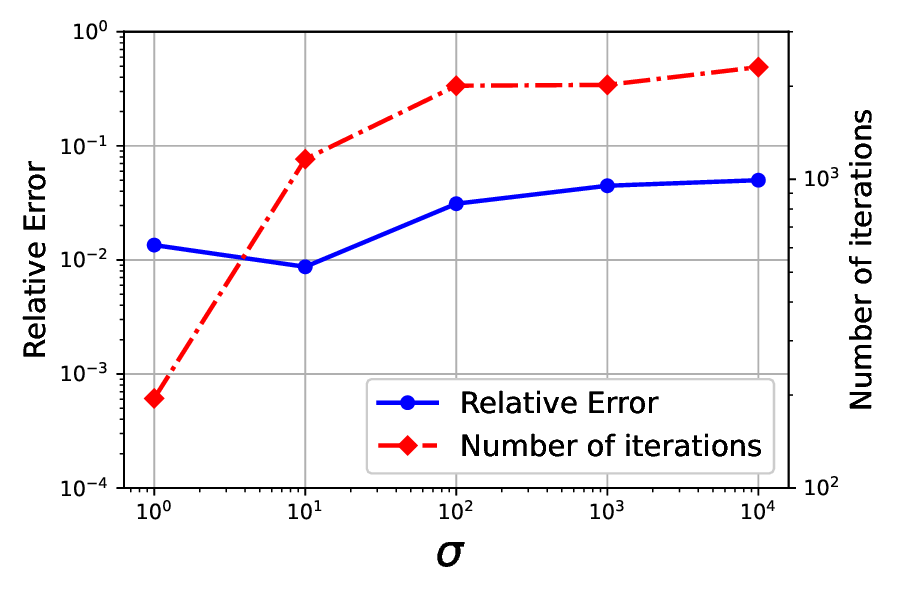}
	\caption{Relative error in the computed electric field and the number of solver iterations as functions of conductivity $\sigma$ of the inner layer of the two layered dielectric cylinder.}
	\label{fig:sigma_convergence_plot}
\end{figure}

The third scenario investigates the impact of varying the dielectric contrast ratio by modeling the inner region of the dielectric cylinder as a conductive medium. The discretization density and other simulation parameters of the first example are maintained. The conductivity of the inner region is varied from 1 to $10^{4}$, and the corresponding relative error in the scattered field as well as the number of iterations required for convergence are reported in Fig.~\ref{fig:sigma_convergence_plot}. While the relative error remains within an acceptable range with the increasing conductivity, the number of iterations increases rapidly up to $\sigma = 10^{2}$, after which the growth becomes more gradual. 


\section{Conclusion}
The proposed reformulated integral equation has been effectively validated through three numerical experiments involving a 2D two-layered dielectric cylinder. This work accurately captures the scattering behavior under TE polarized plane wave illumination as evidenced by the excellent agreement between the computed and analytical solutions of the scattered field. The discretization using RWG basis functions over a triangular mesh proved accurate and effective. Convergence studies further demonstrated that increasing the discretization density leads to improved accuracy of the numerical solution, while the iteration count remains largely unaffected by mesh refinement. This highlights that the proposed discretization preserves the stability and convergence of the underlying integral equation with respect to discretization density.

Finally, the evaluation of the formulation to varying dielectric contrast was examined by modeling the inner region as a conductive medium. Although the relative error remained within acceptable bounds, the number of iterations increased significantly with conductivity. This behavior is consistent with known issues \cite{henry2022low} associated with the underlying integral equations in media with high dielectric contrast and may lead to a computational bottleneck in the forward step of inverse scattering algorithms. The overall results corroborates the effectiveness and applicability of this work for modeling complex inhomogeneous dielectric media.

\bibliographystyle{IEEEtran}
\bibliography{IEEEabrv,Bibliography}


\end{document}